\documentclass[a4paper,12pt]{amsart}
\usepackage{enumerate}
\usepackage{comment}
\usepackage{amsthm,amsmath, amssymb}
\usepackage[top=35truemm, bottom=35truemm, left=30truemm, right=30truemm]{geometry}

\newcommand{\Assouad}{\dim_{\mathrm{A}}}
\newcommand{\Ubox}{\overline{\dim}_{\mathrm{B}}}

\newcommand{\Lbox}{\underline{\dim}_{\mathrm{B}}\,}
\newcommand{\Uzeta}{\overline{\dim}_{\zeta}\,}
\newcommand{\Lzeta}{\underline{\dim}_{\zeta}\,}
\newcommand{\upD}{\overline{D}}
\newcommand{\lowD}{\underline{D}}

\newcommand{\asdim}{\mathrm{asdim}}

\renewcommand{\epsilon}{\varepsilon}
\renewcommand{\limsup}{\varlimsup}
\renewcommand{\liminf}{\varliminf}

\numberwithin{equation}{section}

\newtheorem{theorem}{Theorem}[section]
\newtheorem{lemma}[theorem]{Lemma}
\newtheorem{corollary}[theorem]{Corollary}

\newtheorem{proposition}[theorem]{Proposition}

\newtheorem{question}[theorem]{Question}
\newtheorem{example}[theorem]{Example}
\theoremstyle{definition}

\newtheorem{remark}[theorem]{Remark}
\newtheorem{notation}[theorem]{Notation}

\title{New fractal dimensions and some applications to arithmetic patches}

\author[K. Saito]{ Kota Saito }
\address{Kota Saito\\
Graduate School of Mathematics\\ Nagoya University\\ Furocho\\ Chikusa-ku\\ Nagoya\\ 464-8602\\ Japan }
\curraddr{}
\email{m17013b@math.nagoya-u.ac.jp}
\thanks{}



\date{}

\dedicatory{}

\begin{document}
\maketitle

\begin{abstract} 
In this paper, we define new fractal dimensions of a metric space in order to calculate the roughness of a set on a large scale. These fractal dimensions are called upper zeta dimension and lower zeta dimension.  The upper zeta dimension is an extension of the zeta dimension introduced by Doty, Gu, Lutz, Elvira, Mayordomo, and Moser. We show that the upper zeta dimension is always a lower bound for the Assouad dimension. Moreover, we apply the upper zeta dimension to the existence of weak arithmetic patches of a given set. Arithmetic patches are higher dimensionalized arithmetic progressions. As a corollary,  we get the affirmative solution to a higher dimensional weak analogue of the Erd\H{o}s-Tur\'an conjecture. Here the one dimensional case is proved by Fraser and Yu. As examples, we prove the existence of weak arithmetic patches of the set of all irreducible elements of $\mathbb{Z}[\alpha]$, and the set of all prime numbers of the form $p(p(n))$, where $\alpha$ is an imaginary quadratic integer and $p(n)$ denotes the $n$-th prime number.     
\end{abstract}

\section{Introduction} 
Fractal dimensions are main methods to analyze the roughness of a set. In this paper, we analyze the roughness of a set on a large scale. The first goal of this paper is extending the domain of the \textit{zeta dimension} to a general metric space. The second is comparing the extended zeta dimension to the Assouad dimension. The third is applying to the existence of higher dimensional `weak' arithmetic progressions of given sets. The \textit{zeta dimension} was introduced by Doty, Gu, Lutz, Mayordomo, and Moser \cite{DGLMM}. For every $A\subseteq \mathbb{Z}^d$, the \textit{zeta dimension} of $A$ is defined as
\begin{equation}\label{DotyZetaDim}
	\mathrm{Dim}_{\zeta} A= \inf \{\sigma \geq 0\ : \ \zeta_{A}(\sigma) <\infty  \},
\end{equation}
where the function $\zeta_A(\sigma)$ is defined by
\begin{equation}\label{Zeta}
	\zeta_A(\sigma)=\sum_{v\in A\setminus \{0\}} \|v\|^{-\sigma},
\end{equation}
and $\|\cdot \|$ denotes the Euclidean norm on $\mathbb{R}^d$. Note that the domain of the zeta dimension is only the family of all subsets of $\mathbb{Z}^d$. For every discrete set $F\subset \mathbb{R}^d$, 
\begin{equation}
e(F) = \inf \{\sigma\geq 0 \colon \sum_{v\in F \setminus \{0\}} \| v\|^{-\sigma} <\infty \}
\end{equation}    
is called the \textit{exponent} of $F$ \cite{EssouabriLichtin}. It is clear that $e(F)=\mathrm{Dim}_\zeta F$ if $F$ is a subset of $\mathbb{Z}^d$. Thus $e(F)$ can be considered as an extension of the zeta dimension. The first goal of this paper is to give another extension of the zeta dimension for a general metric space. Let us define new fractal dimensions. Let $(X,d)$ be a metric space, and fix $\alpha \in X$. For every $F\subseteq X$, the \textit{upper zeta dimension} of $F$ is defined by
\begin{align}\label{UzetaDef}
\Uzeta F = \lim_{r\rightarrow \infty} \inf \Big\{\;\sigma\geq 0 & \colon \exists C>0\ \exists R_0>0\ \forall R\geq R_0 
\\ \nonumber
& \, \hspace{2cm} N\Bigl(B(\alpha, R)\cap F,r\Bigr)\leq  C R^\sigma \Big\},
\end{align}

and the \textit{lower zeta function} of $F$ is defined by 
\begin{align}\label{LzetaDef}
\Lzeta F = \lim_{r\rightarrow \infty} \sup \Big\{\;\sigma\geq 0 \ & \colon \exists C>0\ \exists R_0>0\ \forall R\geq R_0 \ 
\\ \nonumber
&\, \hspace{2cm} N \Bigl(B(\alpha, R)\cap F,r\Bigr)\geq  C R^\sigma\Big\}.
\end{align}
Here $B(\alpha, R)$ denotes that the closed ball centered at $\alpha\in X$ with radius $R>0$, and $N(E,r)$ denotes the smallest number of sets with diameter less than or equal to $r$ required to cover $E$. In Section~\ref{well-definedness}, we will show that the limits on the right hand sides of (\ref{UzetaDef}) and (\ref{LzetaDef}) go to infinity or converge, moreover, the upper and lower zeta dimensions do not depend on the choice of any $\alpha\in X$. 

\begin{theorem}\label{main2} For every set $F\subset \mathbb{R}^d$ satisfying $\inf \{\|x-y\|: x,y\in F, x\neq y\}>0$, we have
\[
\Uzeta F=e(F),
\]
In particular, for all $F\subseteq \mathbb{Z}^d$, we have
\[
\Uzeta F =\mathrm{Dim}_\zeta F.
\]
\end{theorem}
The upper zeta dimension can be considered as an extension of the classical zeta dimension (\ref{DotyZetaDim}) from Theorem~\ref{main2}.

The second goal of this paper is to find the inequalities among the upper zeta, lower zeta, and  Assouad dimensions. The Assouad dimension was introduced by Assouad \cite{Assouad}. The precise definition of the Assouad dimension will be given in Section~\ref{ProofofThm}. We refer \cite{Fraser, Robinson} to the reader who wants to know details. For every metric space $X$, $\Assouad X$ denotes the Assouad dimension of $X$. As a result, we get the following relations: 
 \begin{theorem}\label{main4} For all metric space $(X,d)$, we have
 \[
 \Lzeta X \leq \Uzeta X \leq \Assouad X.
 \]
 \end{theorem}
 This theorem is an analogue of the formula
 \begin{equation}\label{boxineq}
 \Lbox X \leq \Ubox X \leq \Assouad X,
 \end{equation}
 where $X$ is any bounded metric space. Here we define 
 \begin{gather}
 \Lbox X = \liminf_{\delta \rightarrow +0} \frac{\log N(X,r)}{-\log r} ,\quad \Ubox =\limsup_{r\rightarrow +0} \frac{\log N(X,r)}{-\log r} 
 \end{gather}
 for every metric space $X$. We say that $\Lbox X$ is the lower box dimension of $X$, and $\Ubox X$ is the upper box dimension of $X$. The formula (\ref{boxineq}) can be seen in \cite{Fraser, Robinson}. The Assouad dimension is always an upper bound for the Hausdorff dimension. We refer the reader to \cite{Falconer, Robinson} for more details on Hausdorff and box dimensions. In Appendix~\ref{appendixA}, we will show that the upper and lower zeta dimensions share some properties with the upper and lower box dimensions, respectively.

\begin{theorem}[quasi-isometric invariance]\label{quasi-iso}
Let $(X,d)$ and $(X',d')$ be metric spaces. If $f:(X,d)\rightarrow (X',d')$ is quasi-isometric \textit{i.e.} there exist $C\geq 1$ and $K\geq 0$ such that for all $x,y\in X$
\[
	\frac{1}{C}d(x,y)-K \leq d'(f(x),f(y)) \leq C d(x,y)+K
\]
holds, then we have
\[
	\Uzeta f(X) = \Uzeta X,\quad \Lzeta f(X) = \Lzeta X.
\]
\end{theorem} 
The asymptotic dimension which is introduced by Gromov \cite{Gromov} is also a quasi-isometrically invariant dimension. Other classical dimensions such as lower, Hausdorff, packing, lower box, upper box, and Assouad dimensions are not invariant with respect to quasi-isometric since a map from a bounded open set to a singleton is quasi-isometric. We will give the precise definition of the asymptotic dimension of a metric space in Appendix~\ref{AsympZeta}. This dimension is considered as an large-scale analogue of the topological dimension. We refer \cite{BellDranishnikov} to the reader for more details. On the other hand, the lower and upper zeta dimensions can be considered as large-scale analogues of the lower and upper box dimensions, respectively, since it follows that
\begin{equation}\label{zetadimsonR}
 \Uzeta F = \limsup_{R\rightarrow \infty} \frac{\log N\big(B(\alpha, R)\cap F,r\big)}{\log R}, \quad
 \Lzeta F = \liminf_{R\rightarrow \infty} \frac{\log N\big(B(\alpha, R)\cap F,r\big)}{\log R}
\end{equation}
for every $F\subseteq \mathbb{R}^d $, $\alpha \in \mathbb{R}^d$, and $r>0$. Note that the zeta dimensions for a subset of $\mathbb{R}^d$ do not depend on the choice of any $\alpha$ or $r$. We will prove the formula (\ref{zetadimsonR}) in Section~\ref{well-definedness}. The topological dimension of a compact set is always a lower bound for the lower box dimension of the set. Thus we might guess that the asymptotic dimension would be always a lower bound for the lower zeta dimension. However, it fails 
in general. In Appendix~\ref{AsympZeta}, we will construct an example of a set of which the asymptotic dimension equals $1$, but the lower zeta dimension equals $0$. We do not find any relations between the asymptotic dimension and the zeta dimensions in this paper. Thus we propose the following remaining question:
\begin{question}
Are there any relations between the asymptotic dimension and zeta dimensions?
\end{question}

The third goal of this paper is to apply the upper zeta dimension to the existence of higher dimensional `weak' arithmetic progressions. Here we now define an arithmetic patch which is a higher dimensionalized arithmetic progression \cite{FraserYu1}. Let $\mathbf{e}=\{e_1,\ldots, e_m\}$ be a linearly independent set of vectors in $\mathbb{R}^d$ where $1\leq m\leq d$. For every $k\in \mathbb{N}$ and $\Delta>0$ we say that a set $P\subset \mathbb{R}^d$ is an \textit{arithmetic patch (AP)} of size $k$ and scale $\Delta$ with respect to orientation $\mathbf{e}$ if 
\[
	P=\left\{t+\Delta\sum_{i=1}^m x_ie_i \ : \ x_i=0,1,\ldots,k-1 \right\}
\] 
for some $t\in \mathbb{R}^d$. For every $\epsilon\in [0,1)$, we say that $Q \subset \mathbb{R}^d$ is a $(k, \varepsilon, \mathbf{e})${\it -AP} if there exists an arithmetic patch $P$ of size $k$, and scale $\Delta >  0$ with respect to orientation $\mathbf{e}$ such that
\begin{equation}\label{f1.1}
\# Q=\# P 
\end{equation}
and
\begin{equation}\label{f1.2}
\sup_{x\in P} \inf_{y \in Q} \|x-y\| \leq \epsilon \Delta. 
\end{equation}
 Note that $(k,0,\{1\})$-APs of real numbers are ordinary arithmetic progressions of length $k$. Fraser and Yu gave the original notion of $(k,\epsilon,\mathbf{e})$-APs in \cite{FraserYu1}. The term `$(k,\epsilon,\mathbf{e})$-APs' was firstly seen in \cite{FraserSaitoYu}. The existence of $(k,\epsilon, \mathbf{e})$-APs of a given set $F$ is strongly connected with the Assouad dimension of $F$. Fraser and Yu showed that a subset of $\mathbb{R}^d$ has Assouad dimension $d$ if and only if the set contains $(k, \varepsilon, \mathbf{e})$-APs for every $k\geq 3$, $\epsilon>0$, and basis $\mathbf{e}$. Note that Fraser and Yu say that $F$ aymptotically contains arbitrarirly large arithemetic patches in \cite{FraserYu1} if $F$ contains $(k, \varepsilon, \mathbf{e})$-APs for every $k\geq 3$, $\epsilon>0$ where $\mathbf{e}$ denotes some fixed basis on $\mathbb{R}^d$. Furthermore, Fraser, the author, and Yu gave the quantitative upper bound of the Assouad dimension of a subset of $\mathbb{R}^d$ which does not contain $(k,\epsilon, \mathbf{e})$-APs. We will restate this result more precisely in Section~\ref{ProofofCoro}. 

 Recently, problems finding arithmetic progressions get interests from many researchers. For examples, Szemer\'edi showed that any subset of positive integers with positive upper density contains arbitrarily long arithmetic progressions, which is called Szemer\'edi's theorem \cite{Szemeredi}.  Furstenberg and Katznelson gave an higher dimensional extension of Szemer\'edi's theorem \cite[Theorem~B]{FurstenbergKatznelson}. Furthermore, Green and Tao showed that the set of all prime numbers contains arbitrarily long arithmetic progressions \cite{GreenTao}, which is called the Green-Tao theorem. Tao also gave an extension of the Green-Tao theorem for the set of all primes in $\mathbb{Z}[i]$. From his result, the set of all primes in $\mathbb{Z}[i]$ contains $(k,0,\{1,i\})$-APs for every $k\geq 2$. One of the biggest conjectures on the existence of arithmetic progressions is the Erd\H{o}s-Tur\'an conjecture. This conjecture states that a subset of positive integers whose sum of reciprocals is divergent would contain arbitrarily long arithmetic progressions. Note that the Erd\H{o}s-Tur\'an conjecture is still open even if the length of arithmetic progressions is equal to 3.
 
 \begin{theorem}[Fraser and Yu \cite{FraserYu1}]\label{ErdosFraserYu}
Any subset of positive integers whose sum of reciprocals is divergent contains $(k,\epsilon, \{1\})$-APs for every $k\geq 3$ and $\epsilon>0$.
\end{theorem}
 
  Theorem~\ref{ErdosFraserYu} can be considered as a weak affirmative solution to the Erd\H{o}s-Tur\'an conjecture. Fraser gave an elementary proof of Theorem~\ref{ErdosFraserYu} in \cite{Fraser2}. If we could replace $\epsilon>0$ to $\epsilon=0$ in Theorem~\ref{ErdosFraserYu}, we got the affirmative solution to the Erd\H{o}s-Tur\'an conjecture. However, we do not get results on the case when $\epsilon=0$ in this paper. As an application, we get the following higher dimensional and quantitative extension of Theorem~\ref{ErdosFraserYu};
\begin{corollary}\label{FraserSaitoYuET} Let $F$ be a subset of $\mathbb{R}^d$ satisfying $\inf \{\|x-y\|: x,y\in F, x\neq y\}>0$.  for all distinct elements $x,y\in F$. Fix $k\geq 3$, $\epsilon>0$, and $1\leq m\leq d$. If the series 
\[
\sum_{v\in F\setminus \{0\}} \|v\|^{-\sigma_0} 
\]
is divergent for some
\[
\sigma_0 > d+\frac{\log (1-1/k^m)}{\log (k \lceil \sqrt{d}/(2\epsilon) \rceil) },
\]
then $F$ contains $(k,\epsilon, \mathbf{e})$-APs for every a set of orthogonal unit vectors $\mathbf{e}=\{e_1,e_2,\ldots , e_m\}$. 
\end{corollary}
Moreover, we can get the following corollary:
\begin{corollary}\label{main5}
Let $F$ be a subset of $\mathbb{Z}^d$ such that
\begin{equation}\label{introF1}
\sum_{v\in F\setminus \{0\}} \|v\|^{-\sigma}=\infty
\end{equation}
holds for every $0<\sigma<d$. Then $F$ contains $(k,\epsilon, \mathbf{e})$-APs for every $k\geq 3$,   $0<\epsilon <1$, and orthonormal basis $\mathbf{e}$.
\end{corollary} 
 Fraser and Yu have already given another higher dimensional extension of Theorem~\ref{ErdosFraserYu} \cite{FraserYu2}. They replace $\epsilon \Delta$ to $\Delta^\alpha$ in (\ref{f1.2}) where $\alpha\in (0,1)$. Neither the extention by Fraser and Yu implies Corollary~\ref{main5}, nor does Corollary~\ref{main5} imply the extension by Fraser and Yu. 
 
 As applications of Corollary~\ref{main5}, in Section~\ref{Examples} we will show the existence of $(k,\epsilon, \mathbf{e})$-APs of the set of all irreducible elements of $\mathbb{Z}[\alpha]$, and the set of all prime numbers of the form $p(p(n))$, where $\alpha$ is an imaginary quadratic integer and $p(n)$ denotes the $n$-th prime number.

\begin{question}\label{higherdimET}
Let $\mathbf{e}$ be the standard basis of $\mathbb{R}^d$. 
Is it true that any subset $F$ of $\mathbb{Z}^d$ satisfying
\[
\sum_{v\in F\setminus\{0\} } \|v\|^{-d}=\infty
\]
contains $(k,0,\mathbf{e})$-APs for every $k\geq 2$?

\end{question}
This is a higher dimensional extension of the Erd\H{o}s-Tur\'an conjecture. The answer to Question~\ref{higherdimET} would be `yes' in the view of Corollary~\ref{main5}. In particular, this was conjectured by Graham when $d=2$ \cite[Conjecture~11.5.7]{Graham}.  

\begin{notation} Let $(X,d)$ be a metric space. We give the following notations:
\begin{itemize}
\item For every set $A\subseteq [0,\infty)$, $\inf A$ denotes the maximum of lower bounds of $A$ where we define $\inf A=\infty$ if $A=\emptyset$, and $\sup A$ denotes the minimum of upper bounds of $A$ where we define $\sup A=0$ if $A=\emptyset$ and  $\sup A=\infty$ if $A$ does not have any upper bounds; 
\item The diamiter of $U\subseteq X$ is defined by $\sup\{d(x,y) \colon x,y\in U \}$;
\item For every set $A$, $|A|$ denotes the cardinality of $A$;
\item Let $E$ be a subspace of $X$. A finite family $\mathcal{U}$ of sets is called an $r$-cover of $E$ if the union of all sets in $\mathcal U$ contains $E$ and the diameter of $U$ is less than or equal to $r$ for every $U\in \mathcal{U}$; 
\item We can rewrite $N(E,r)=\min\{|\mathcal{U}|\colon \mathcal{U} \text{ is an $r$-cover of $E$} \} $, where $N(E,r)=\infty$ if there does not exist any $r$-covers of $E$.
\end{itemize}
\end{notation}

\section{Well-definedness of the upper and lower zeta dimensions}\label{well-definedness}

Let $(X,d)$ be a metric space. For every $F\subseteq X $, $\alpha\in X$, and $r>0$,   we define
 \begin{equation}\label{upDDef}
  \upD(F;\alpha, r)= \inf \Big\{\;\sigma\geq 0 \colon \exists C>0\ \exists R_0>0\ \forall R\geq R_0 \ 
N\Bigl(B(\alpha, R)\cap F,r\Bigr)\leq  C R^\sigma \Big\},
  \end{equation}
 and
 \begin{equation}\label{lowDDef}
\lowD(F; \alpha, r)= \sup \Big\{\;\sigma\geq 0 \  \colon \exists C>0\ \exists R_0>0\ \forall R\geq R_0 \ 
N \Bigl(B(\alpha, R)\cap F,r\Bigr)\geq  C R^\sigma\Big\}.
 \end{equation}
 We firstly show the well-definedness of the upper and lower zeta dimensions \textit{i.e.} the limits 
 \begin{gather*}
 \lim_{r\rightarrow \infty} \upD(F;\alpha, r),\ \quad \lim_{r\rightarrow \infty}	\lowD(F; \alpha, r)
 \end{gather*}
  exist or go to infinity. Note that the limits of $\upD(F;\alpha, r)$ and $\lowD(F; \alpha, r)$ with respect $r$ are formally equal to $\Uzeta F$ and $\Lzeta F$, respectively.
Before proof, we show the following lemma:
  
\begin{lemma}\label{equivforms1}Let $(X,d)$ be a metric space. For every $F\subseteq X $, $\alpha\in X$, and $r>0$, we have
\begin{equation}\label{equivup1} 
 \upD (F;\alpha,r) = \limsup_{R\rightarrow \infty} \frac{\log N\big(B(\alpha, R)\cap F,r\big)}{\log R}
\end{equation}
and
\begin{equation}\label{equivlow1}
 \lowD (F; \alpha, r) = \liminf_{R\rightarrow \infty} \frac{\log N\big(B(\alpha, R)\cap F,r\big)}{\log R}.
\end{equation}

 \end{lemma}
 
 \begin{proof}
  Fix $r>0$ and let $\sigma_0=\upD (F; \alpha, r)$. If $\sigma_0<\infty$, then for any fixed $\varepsilon>0$ there exist 
  $ C=C(r,\epsilon)>0$ and  $ R_0=R_0(r,\epsilon)>0$ such that 
$
   N(B(\alpha, R)\cap F,r )\leq  CR^{\sigma_0+\varepsilon}
$
holds for all $R\geq R_0$. Thus by taking the logarithm on both sides,  we obtain
\[
	\frac{\log{N(B(\alpha,R)\cap F,r)}}{\log R}\leq \sigma_0 + \varepsilon + \frac{C_{r}}
	{\log{R}} 
\]  
for all $R\geq R_0$. Therefore we get
\[
	\limsup_{R\rightarrow \infty} \frac{\log N (B(\alpha, R)\cap F,r)}{\log R}\leq 
	\sigma_0\;.
\]
Conversely, let
\[
	\sigma'_0=\limsup_{R\rightarrow \infty} \frac{\log N(B(\alpha, R)\cap F,r)}{\log 
	R}.
\]
If $\sigma'_0<\infty$, then for every $\varepsilon>0$ there exists $R_0>0$ such that
\[
	\frac{\log N(B(\alpha, R)\cap F,r)}{\log R}\leq \sigma_0'+\varepsilon
\]
for all $R\geq R_0$. Thus it follows that $N(B(\alpha, R)\cap F,r)\leq R^{\sigma_0+\varepsilon}$, which implies $\upD_\alpha(F,r)\leq \sigma'_0+\varepsilon$ from (\ref{upDDef}). By the above proof, $\sigma_0=\infty$ is equivalent to $\sigma'_0=\infty$.  Therefore we have (\ref{equivup1}). Similarly we have (\ref{equivlow1}). 
\end{proof}
 
 Thanks to Lemma~\ref{equivforms1}, we can easily see that 
 \begin{gather}\label{monoD}
 \upD(F;\alpha, \ell) \geq \upD(F; \alpha, r),\quad \lowD(F;\alpha, \ell) \geq \lowD(F;\alpha,r)
 \end{gather}
 for every $0<\ell \leq r$ since it follows that
 \begin{equation}\label{f1section3}
 N (B(\alpha, R)\cap F,\ell )\geq  N (B(\alpha, R)\cap F,r)
 \end{equation}
 from the definition. By combining (\ref{f1section3}) and Lemma~\ref{equivforms1}, we obtain the inequalities (\ref{monoD}). Furthermore, the inequalities (\ref{monoD}) imply that
  \begin{gather*}
 \lim_{r\rightarrow \infty} \upD(F;\alpha, r)=\Uzeta F, \quad \lim_{r\rightarrow \infty} \lowD(F; \alpha, r)=\Lzeta F
 \end{gather*}
  exist or go to infinity. Hence the upper and lower zeta dimensions are well-defined.

 We next see that the upper and lower zeta dimensions do not depend on the choice of any $\alpha\in X$, that is, the following lemma holds:

\begin{lemma}\label{indalpha}
 Let $(X,d)$ be a metric space. We have
\begin{gather}
 \upD(F; \alpha, r)=\upD(F; \beta, r),\quad \lowD(F; \alpha, r)=\lowD(F;\beta, r)
\end{gather}
for every $F\subseteq X$, $\alpha, \beta\in X$, and $r>0$.
\end{lemma}

\begin{proof}
Let $l=d(\alpha,\beta)$. It follows that $B(\beta,R)\subset B(\alpha,R+l)$ for every $R>0$ since if $x \in B(\beta, R)$ holds, then we obtain 
$
	d(\alpha,x)\leq d(\alpha,\beta)+d(\beta,x) \leq l+ R
$ by the triangle inequality. Hence we have
\[
	N\big(B(\beta,R)\cap F, r\big) \leq N\big(B(\alpha, R+l)\cap F,r\big)\;
\]
by definition. By Lemma~\ref{equivforms1}, this yields that
\begin{align*}
\upD(F; \beta, r)&=\limsup_{R\rightarrow \infty} \frac{\log N(B(\beta, R)\cap F,r)}{\log R}\\
&\leq \limsup_{R\rightarrow \infty} \frac{\log N(B(\alpha, R+l)\cap F,r )}{\log (R+l)}\cdot \frac{\log{(R+l)}}{\log R} =\upD(F; \alpha,r).
\end{align*}
Similarly, we have the opposite inequality $\upD(F;\beta , r)\geq\upD(F;\alpha,r)$. Therefore $\upD(F;\beta,r)=\upD(F,\alpha, r)$ holds. By the same way, we obtain $\lowD(F; \alpha,r)=\lowD(F; \beta, r)$. 
\end{proof}

 Furthermore, we will show that $\upD (F;\alpha,r)$ and $\lowD(F;\alpha,r)$ do not depend on the choice of any $r>0$ if the Assouad dimension of $X$ is finite. Here the \textit{Assouad dimension} of $F\subseteq X$ is defined by
\begin{align}
	\Assouad F = \inf \Big\{\;\sigma\geq 0 &\colon \exists C>0\ \forall r>0\ \forall R> r \ \forall x\in F\\ \nonumber
	&\, \hspace{2cm}\  N\Bigl(B(\alpha, R)\cap F,r\Bigr)\leq  C R^\sigma \Big\}.
\end{align}
\begin{lemma}\label{indr}Let $(X,d)$ be a metric space of which the Assouad dimension is finite. Fix $\alpha \in X$ We have 
\begin{gather}
	\upD (F;\alpha, \ell)=\upD (F;\alpha, r)\text{ and } \lowD(F;\alpha, \ell)=\lowD(F;\alpha, r)
\end{gather}
for every $\ell, r>0$, and $F\subseteq X$.
\end{lemma} 

\begin{proof}  Fix $F\subseteq X$. Without loss of generality, we may assume that $0< \ell \leq r$. Thus by the inequality (\ref{monoD}), it is enough to  show that
\begin{gather}\label{f3section3}
\upD ( F;\alpha,\ell) \leq \upD ( F;\alpha,r) \quad\text{and}\quad \lowD ( F;\alpha,\ell) \leq \lowD ( F;\alpha,r)
\end{gather}
for every $0<\ell\leq r$. Fix real numbers $\ell$ and $r$ with $0<\ell\leq r$, and fix $F\subseteq X$. From the condition $\Assouad X<\infty$ and the monotonicity of the Assouad dimension, it follows that $\Assouad F <\infty$. Thus we can take an $r$-cover $\{U_j\}_{j=1}^N$ of $B(\alpha, R)\cap F$, where $N=N(B(\alpha, R)\cap F,r)$. Each $U_j$ intersects with $B(\alpha,R)\cap F$ because of the minimality of $N$. From the triangle inequality, there exists a family $\{B_j\}_{j=1}^N$ of closed balls centered at an element in $F$ with radius $r$ such that $U_j\subseteq B_j$ for all $j=1,\ldots, N$. Fix $j=1,\ldots, N$. From the finiteness of the Assouad dimension of $F$, there exists $C>0$ such that 
\[
N(B_j\cap F, \ell)\leq C(r/\ell)^{\Assouad F+1},  
\]
where $C$ does not depend on $j, r,$ or $\ell$. Hence $B(\alpha ,R)\cap F$ can be covered by at most $C(r/\ell)^{\Assouad F+1}N$ sets with diameter less than or equal to $\ell$. Hence we obtain
\[
N(B(\alpha, R)\cap F, \ell) \leq  C(r/\ell)^{\Assouad F+1}N= C(r/\ell)^{\Assouad F+1}N(B(\alpha, R)\cap F,r).
\]  
By Lemma~\ref{equivforms1}, we conclude (\ref{f3section3}).
\end{proof}

\begin{remark}\label{r-ind} Let $(X,d)$ be a metric space of which the Assouad dimension is finite. Thanks to Lemma~\ref{indalpha} and Lemma~\ref{indr}, we obtain that  
\begin{equation*}
\upD(F; \alpha, r)=\Uzeta F, \quad \lowD(F;\alpha, r)=\Lzeta F
\end{equation*}
for every $F\subseteq X$, $\alpha\in X$, and $r>0$. In particular, if $(X,d)$ is a $d$-dimensional Banach space, then the Assouad dimension of $X$ is equal to $d$. Therefore by Lemma~\ref{equivforms1} we get the equivalent forms
\[ 
 \Uzeta F = \limsup_{R\rightarrow \infty} \frac{\log N\big(B(\alpha, R)\cap F,r\big)}{\log R}, \quad
 \Lzeta F = \liminf_{R\rightarrow \infty} \frac{\log N\big(B(\alpha, R)\cap F,r\big)}{\log R}
\]
for every $F\subseteq X $, $\alpha \in X$, and $r>0$. Furthermore, if $F\subseteq X$ satisfies that there exists $r>0$ such that $d(x,y)>r$ for all $x,y\in F$ with $x\neq y$, then we have
\[
N(B(\alpha, R)\cap F, r/2)=|B(\alpha,R)\cap F| 
\]
for every $\alpha \in X$, which implies that
\begin{equation}\label{Euclidsp}
 \Uzeta F = \limsup_{R\rightarrow \infty} \frac{\log |B(\alpha, R)\cap F|}{\log R}, \quad
 \Lzeta F = \liminf_{R\rightarrow \infty} \frac{\log |B(\alpha, R)\cap F|}{\log R}
\end{equation}
for every $\alpha\in X$.
\end{remark}

\section{Proof of Main Results}\label{ProofofThm}

\begin{proof}[Proof of Theorem~\ref{main4}] 
By Lemma~\ref{equivforms1}, then it is clear that $\Lzeta X\leq \Uzeta X$. We show that $\Uzeta X\leq \Assouad X$. If $\dim_\mathrm{A} X=\infty$ holds, it is clear that $\Uzeta X \leq \dim_\mathrm{A} X$. Thus we may assume that $\Assouad X<\infty$. We choose any $\sigma$  such that $\sigma>\Assouad X$.  From the definition of the Assouad dimension, there exists $C>0$ such that for all $r>0$, $R>r$, and $x\in X$ we have 
\[
N(B(x,R)\cap X ,r) \leq C\left(\frac{R}{r}\right)^\sigma
\]

We now fix $r>0$ and $\alpha\in X$. Then the following inequality holds:
\[
N(B(\alpha,R)\cap X,r)\leq   \ C \left(\frac{R}{r}\right)^\sigma \leq Cr^{-\sigma} R^{\sigma}.
\]
By Lemma~\ref{equivforms1}, we conclude $\Lzeta X \leq \Uzeta X\leq \Assouad X$.
\end{proof}

\begin{proof}[Proof of Theorem~\ref{quasi-iso}]
We write $B'(\alpha, R)$ as the closed ball of $X'$ centered at $\alpha \in X'$ with radius $R$. Fix any sufficiently large $R>0$ and $r>0$. We take any $(r/C-K)$-cover $\{V_{j}\}_{j=1}^N$ of $B'(f(\alpha), CR+K)\cap f(X)$, where let $N=N(B'(f(\alpha), CR+K), r/C-K)$. Then
\[
	\{f^{-1}(V_1),\ldots, f^{-1}(V_N)\}
\]
can cover $B(\alpha, R)\cap X$ since it follows that
\[
\bigcup_{j=1}^N f^{-1}(V_j) \supseteq f^{-1}\left(\bigcup_{j=1}^N V_j\right)\supseteq f^{-1}(B'(f(\alpha), CR+K)\cap f(X)) \supseteq B(\alpha, R)\cap X.
\]
Furthermore, the diameter of $f^{-1}(V_j)$ is less than or equal to $r$ for every $j=1,2,\ldots, N$. Hence we have
\begin{equation}\label{quasi-iso-f1}
	N(B(\alpha, R)\cap X, r) \leq N(B'(f(\alpha), CR+K)\cap f(X), r/C-K). 
\end{equation}
Therefore Lemma~\ref{equivforms1} and (\ref{quasi-iso-f1}) imply that
\[
	\Uzeta f(X) \geq \Uzeta X,\quad \Lzeta f(X) \geq \Lzeta X.
\] 
We next take any $((r-K)/C)$-cover $\{V_{j}\}_{j=1}^N$ of $B(\alpha, C(R+K))\cap X$, where let $N=N(B(\alpha,C(R+K))\cap X, (r-K)/C)$. Then
\[
	\{f(V_1),\ldots, f(V_N)\}
\]
can cover $B'(f(\alpha), R)\cap f(X)$ since it follows that

\[
\bigcup_{j=1}^N f(V_j) \supseteq f\left(\bigcup_{j=1}^N V_j\right)\supseteq f(B(\alpha,C(R+K))\cap X) \supseteq B'(f(\alpha), R)\cap f(X).
\]
The diameter of $f(V_j)$ is less than or equal to $r$ for every $j=1,2,\ldots, N$. Hence we have
\begin{equation}\label{quasi-iso-f2}
	N(B'(f(\alpha), R)\cap f(X), r) \leq N(B(f(\alpha), C(R+K))\cap f(X), (r-K)/C). 
\end{equation}
Lemma~\ref{equivforms1} and (\ref{quasi-iso-f2}) imply that
\[
	\Uzeta f(X) \leq \Uzeta X,\quad \Lzeta f(X) \leq \Lzeta X.
\]
\end{proof}

Before proving Theorem~\ref{main2}, we define the function $\zeta(\sigma;\alpha, \Gamma)$ and the notion $r$-\textit{net}. We say that a metric space $(\Gamma, d)$ is $r$-discrete if $d(x,y)\geq r$ holds for all $x,y\in \Gamma$ , and we say that $(\Gamma, d)$ is \textit{partially finite} if $B(x,R)\cap \Gamma$ is  finite for all $x\in \Gamma$ and $R>0$. Let $(X,d)$ be a metric space, $(\Gamma,d)$ be a partially finite and $r$-discrete subspace of $(X,d)$ for some $r>0$, and fix any $\alpha\in X$. The function $\zeta(\sigma ;\alpha, \Gamma )$ is defined by
\begin{equation}\label{MetricZeta}
	\zeta (\sigma; \alpha, \Gamma)=\sum_{\gamma\in \Gamma\setminus \{\alpha\}} d(\gamma, \alpha)^{-\sigma}
\end{equation}
for every $\sigma\geq 0$. Note that $\zeta(\sigma; 0, \Gamma)=\zeta_\Gamma(\sigma)$ holds if the metric function $d$ is the Euclidean metric and $\Gamma\subseteq \mathbb{Z}^d$, where the definition of $\zeta_\Gamma(\sigma)$ is given in \ref{Zeta}. Hence the function
 $\zeta (\sigma; \alpha, \Gamma)$ can be considered as an extension of the function $\zeta_\Gamma(\sigma)$ to a partially finite and $r$-discrete metric space. We next say that a  set $\Gamma \subseteq X$ is an $r$-\textit{net} of $X$ for some fixed $r>0$ if we have
\begin{itemize}
\item[(i)]$\Gamma$ is $r$-discrete, and
\item[(ii)] for all $x\in X$ there exists $\gamma\in \Gamma$ such that  $d(x,\gamma)< r$. 
\end{itemize}
From the Zorn's lemma, for every $r>0$, there exists a set $\Gamma\subseteq X$ satisfying (i) and (ii). We will give a proof of the existence of an $r$-net in Appendix~\ref{appendixA}. We say that $\Gamma$ is a \textit{net} of $X$ if $\Gamma$ is an $r$-net of $X$ for some $r>0$. The abscissa of convergence of the function $\zeta (\sigma; \alpha, \Gamma)$  has strong connection to the upper zeta dimension as follows:

\begin{proposition}\label{main1} Let $(X,d)$ be a non-empty metric space of which the Assouad dimension is finite. Fix a net $\Gamma$ of $X$. Then we have
\begin{equation}
\Uzeta X = \inf \{\sigma\geq 0\ \colon\ \zeta(\sigma; \alpha, \Gamma) <\infty \}.
\end{equation}
\end{proposition}

\begin{proof}[Proof of Theorem~\ref{main2} assuming Proposition~\ref{main1}]
Fix an $r$-discrete subset $F$ of $\mathbb{R}^d$ for some $r>0$. Then it is clear that $F$ is a net of $F$. Therefore we have
\[
\Uzeta F = \inf \{\sigma\geq 0\ \colon \zeta(\sigma;0, F) <\infty \}
=\inf \{\sigma\geq 0\ \colon \sum_{v\in F\setminus \{0\} }\|v\|^{-\sigma} <\infty  \} =e(F)
\]
by Proposition~\ref{main1}. In particular, for every $F\subseteq \mathbb{Z}^d$, $e(F)=\mathrm{Dim}_\zeta F$ holds from the definitions of $e(F)$ and $\mathrm{Dim}_\zeta (F)$.
\end{proof}

The remaining part of proving Theorem~\ref{main2} is to show Proposition~\ref{main1}.  

\begin{proof}[Proof of Proposition~\ref{main1}]
Fix $\alpha \in X$. We choose a net $\Gamma$ of $X$. By definition, there exists $r>0$ such that $d(x,y)\geq r$ for all distinct elements $x,y\in \Gamma$, and for all $x\in X$ there exists $\gamma\in \Gamma$ such that $d(x,\gamma)<r$. We firstly show that 
\begin{gather}
\Uzeta \Gamma = \Uzeta X.
\end{gather}
By definition, we see that
\[
	X= \bigcup_{\gamma\in \Gamma} B(\gamma, r)  
\]
By the axiom of choice, there exists a surjective map $f\colon X \rightarrow \Gamma$ such that $x\in B(f(x), r)$ for all $x\in X$.
Then we have
\[
d(x,y)-2r\leq d(f(x),f(y))\leq d(x,y)+2r
\]
by the triangle inequality. Therefore we obtain
\[
\Uzeta \Gamma = \Uzeta X
\]
by Theorem~\ref{quasi-iso}. From Remark~\ref{r-ind} and the condition $\Assouad X<\infty$, it is enough to show that 
\[
\upD(\Gamma ; \alpha, r/2) = \inf\{\sigma\geq 0 \colon \zeta(\sigma ; \alpha, \Gamma)<\infty \}.
\] 
We have $N(B(\alpha, R)\cap \Gamma, r/2)=|B(\alpha,R)\cap \Gamma |$ since $d(x,y)\geq r$ for all distinct elements $x,y\in \Gamma$. Thus by Lemma~\ref{equivforms1}, we have 
\begin{equation}\label{Pr-f1}
\upD(\Gamma;\alpha, r/2)= \limsup_{R\rightarrow \infty}\frac{\log |B(\alpha, R)\cap \Gamma|}{\log R}.
\end{equation}
Let $\sigma_0$ be the right hand side of (\ref{Pr-f1}). It is enough to show that 
\[
\sigma_0=\inf \{ \sigma\geq 0 \:\colon\: \zeta(\sigma; \alpha, F)<\infty \}.   
\]
The condition that $\Assouad X<\infty$ implies $\Assouad \Gamma<\infty$. Therefore $\Gamma$ is partially finite by the definition of the Assouad dimension. Hence we can write
\[
\{d(\alpha, \gamma)\ \colon \ v\in \Gamma\setminus \{\alpha\} \}=\{\ell_j \colon j=0,1,\ldots\}
\]
for some infinite sequence $\ell_0< \ell_1< \cdots$. Let
\begin{gather*}
 	r(j)=|\{ v\in \Gamma \colon \ell_j=d(\alpha, v) \}|, \quad M(\ell) =\sum_{j\: \colon\: \ell_j\leq \ell }r(j).
\end{gather*}
It is clear that $r(j), M(\ell)<\infty$ because $\Gamma$ is partially finite. We fix $\ell>\ell_0$, and choose a non-negative integer $m$ such that $\ell_{m}\leq \ell< \ell_{m+1}$. Then for every $\sigma>0$, we obtain that
\begin{align*}
	&\, \sum_{\substack{v\in F\setminus \{\alpha\}\\   d(\alpha, v)\leq \ell }} \frac{1}{d(\alpha, v)^\sigma}
	 = \sum_{j=0 }^m \frac{r(j)}{\ell_j^\sigma} = \frac{M(\ell_m)}{\ell_m^\sigma}+\sum_{j=0}^{m-1} \left(\frac{1}{\ell_j^\sigma}-\frac{1}{\ell_{j+1}^\sigma}\right)M(\ell_j)\\
	 &= \frac{M(\ell_m)}{\ell_m^\sigma}+\sum_{j=0}^{m-1}  \sigma M(\ell_j)\int_{\ell_j}^{\ell_{j+1}} \frac{1}{t^{\sigma+1}}\: \mathrm{d}t 
	 = \frac{M(\ell_m)}{\ell_m^\sigma}  +\sigma \sum_{j=0}^{m-1}  \int_{\ell_j}^{\ell_{j+1}} \frac{M(t)}{t^{\sigma+1}}\: \mathrm{d}t \\
	 &=  \frac{M(\ell_m)}{\ell_m^\sigma}  +\sigma \int_{\ell_0}^{\ell_{m}} \frac{M(t)}{t^{\sigma+1}}\: \mathrm{d}t.
\end{align*}
Note that this calculation is same as the proof of Abel's identity in \cite[Theorem~4.2]{Apostol}. By the choice of $\sigma_0$, for every $\epsilon>0$ there exists a positive real number $C_\epsilon$ such that 
\[
	M(\ell)\leq  |B(\alpha, \ell)\cap \Gamma| \leq C_\epsilon \ell^{\sigma_0+\epsilon}
\]
holds for all $\ell>0$. Hence we find that
\[
\sum_{\substack{v\in F\setminus \{\alpha\} \\ d(\alpha, v)\leq \ell }} \frac{1}{d(\alpha, F)^{\sigma_0+2\epsilon}} 
= \frac{M(m)}{{\ell_m}^{\sigma_0+2\epsilon}}+\sigma \int_{\ell_0}^{\ell_{m}} \frac{M(\ell)}{\ell^{\sigma_0+2\epsilon+1}}\: \mathrm{d}\ell \leq C_\epsilon \left(  1+\frac{\sigma_0+2\epsilon}{\epsilon} \frac{1}{\ell_0^{\epsilon}} \right),    
\]
which implies that
\[
\sigma_0+2\epsilon \geq \inf \{ \sigma\geq 0 \:\colon\: \zeta(\sigma; \alpha, F) <\infty \}. 
\]
By taking $\epsilon\rightarrow +0$, we obtain
\[
	\sigma_0 \geq \inf \{ \sigma\geq 0 \:\colon\: \zeta(\sigma; \alpha, F)<\infty\}.
\]
We next show the inequality in the other direction.  We take a positive real number $\sigma_1$ such that $\zeta(\sigma;\alpha, \Gamma)$ is convergent. Then let $T=\zeta(\sigma_1;\alpha, \Gamma)$. It is clear that $T\geq M(\ell_j)/\ell_j^{\sigma_1}$ for every non-negative integer $j$ by estimating the partial sum. Fix a real number $\ell>\ell_0$, and choose a non-negative integer $m$ such that $\ell_{m}\leq \ell< \ell_{m+1}$. Then we have
\[
\frac{\log |B(\alpha, \ell)\cap F|}{\log \ell} \leq \frac{\log |B(\alpha, \ell_m)\cap F|}{\log \ell_m} \leq \frac{\log T\ell_m^{\sigma_1}}{\log \ell_m},
\]
which implies that 
\[
\sigma_0=\limsup_{R\rightarrow \infty}\frac{\log |B(\alpha, R)\cap F|}{\log R}\leq \sigma_1.
\]
Therefore we have
\[
 \sigma_0\leq \inf \{ \sigma\geq 0 \:\colon\: \zeta(\sigma; \alpha, F)<\infty\}.
\]
\end{proof}

\section{Proof of Corollary~\ref{FraserSaitoYuET} and Corollary~\ref{main5}} \label{ProofofCoro}

In order to prove Corollary~\ref{FraserSaitoYuET}, we use a quantitative upper bound of the Assouad dimension of a set which does not contain $(k,\epsilon, \mathbf{e})$-APs. This upper bound was given by Fraser, the author, and Yu \cite{FraserSaitoYu}.
\begin{theorem}[Fraser, the author, and Yu \cite{FraserSaitoYu}]\label{FraserSaitoYu}
Fix $k\geq 3$ and $\epsilon>0$. If $F\subseteq \mathbb{R}^d$ does not contain any $(k,\epsilon,\mathbf{e})$-APs for some a set of orthonormal vectors $\mathbf{e}=\{e_1,\ldots, e_m\}$, then
\[
\Assouad F \leq d+\frac{\log (1-1/k^m)}{\log k \lceil \sqrt{d}/(2\epsilon) \rceil }
\]
holds.
\end{theorem} 

\begin{proof}[Proof of Corollary~\ref{FraserSaitoYuET}]
Fix $\alpha =0$.  By the condition, we can take a real number $\sigma_0$ such that
\[
d+\frac{\log (1-1/k^m)}{\log k \lceil \sqrt{d}/(2\epsilon) \rceil }<\sigma_0
\]
and $\zeta(\sigma_0; \alpha, F)=\sum_{v\in F\setminus \{0\}} \|v\|^{-\sigma_0}=\infty$ holds. Therefore we have
\begin{equation}\label{f1section2}
d+\frac{\log (1-1/k^m)}{\log k \lceil \sqrt{d}/(2\epsilon) \rceil }<\sigma_0\leq  \inf\{\sigma \geq 0 \colon \zeta(\sigma; 0, F )<\infty \} 
\end{equation}
By Theorem~\ref{main2} and Theorem~\ref{main4}, we have
\begin{equation}\label{f2section2}
  \inf\{\sigma \geq 0 \colon \zeta(\sigma; 0, F )<\infty \} =\Uzeta F \leq  \Assouad F.
\end{equation}
By combining the inequalities (\ref{f1section2}) and (\ref{f2section2}), we obtain 
\[
d+\frac{\log (1-1/k^m)}{\log k \lceil \sqrt{d}/(2\epsilon ) \rceil }< \Assouad F.
\]
From Theorem~\ref{FraserSaitoYu}, $F$ contains $(k,\epsilon, \mathbf{e})$-APs for every set of orthonormal vectors $\{e_1, \ldots, e_m\}$. 
\end{proof}

\begin{proof}[Proof of Corollary~\ref{main5}]
Fix any $k\geq 3$, $\epsilon>0$, and $d\geq 1$. We substitute $m$ for $d$ in Corollary~\ref{FraserSaitoYuET}.  We observe that
\[
d+\frac{\log (1-1/k^d)}{\log k \lceil \sqrt{d}/(2\epsilon) \rceil }< d
\]
Therefore there exists a real number $\sigma_0$ such that
\[
d+\frac{\log (1-1/k)}{\log k \lceil \sqrt{1}/(2\epsilon) \rceil }<\sigma_0<d
\]
and 
\[
\sum_{n\in F} n^{-\sigma_0}=\infty
\]
 by the condition of $F$. This yields that $F$ contains $(k,\epsilon, \mathbf{e})$-APs for every set of orthonormal vectors $\{e_1, \ldots, e_d\}$ from Corollary~\ref{FraserSaitoYuET}.  
\end{proof}

\begin{remark}
We do not have to restrict a basis $\{e_1,\ldots, e_d\}$ to be orthonormal in Corollary~\ref{main5}. In fact, let $F\subseteq \mathbb{R}^d$. Thanks to the result of Fraser and Yu \cite{FraserYu1}, it follows that $\Assouad F=d$ if and only if $F$ contains $(k,\epsilon, \mathbf{e})$-APs for every $k\geq 3$, $\epsilon>0$, and basis $\mathbf{e}$ of $\mathbb{R}^d$. If $F$ satisfies that
\[
\sum_{v\in F} \|v\|^{-\sigma}=\infty
\]
for every $\sigma<d$ and $F$ is $r$-discrete for some $r>0$, then we have $\Assouad F=d$ by Theorem~\ref{main2} and Theorem~\ref{main4}. Therefore $F$ contains $(k,\epsilon, \mathbf{e})$-APs for every $k\geq 3$, $\epsilon>0$, and basis $\mathbf{e}$ of $\mathbb{R}^d$.
\end{remark}

\section{Examples}\label{Examples}

\begin{lemma}\label{union}
Let $(X,d)$ be a metric space, and let $F$ and $E$  be subspaces of $(X,d)$. Then the following properties hold:
\begin{itemize}
\item[(i)] If $F\subseteq E$, then we have $\Uzeta F\leq \Uzeta E$ and $\Lzeta F \leq \Lzeta E$;
\item[(ii)] $\Uzeta (F\cup E) = \max\{\Uzeta F,\: \Uzeta E \}$.
\end{itemize}
\end{lemma}

\begin{proof}
We first show (i). Fix $\alpha\in X$ and $r>0$. We find that $N(B(\alpha, R)\cap F, r) \leq N(B(\alpha, R)\cap E, r)$ for all $R>0$ by definition.   
Therefore we have $\Uzeta F\leq \Uzeta E$ and $\Lzeta F\leq \Lzeta E$ by Lemma~\ref{equivforms1}. Let us show (ii). It is clear from (i) that
\[
 \Uzeta (F\cup E) \geq \max\{\Uzeta F,\; \Uzeta E \}.
\]
The opposite inequality also holds by the inequality
\[ 
 N(B(\alpha, R)\cap (F\cup E), r) \leq N(B(\alpha, R)\cap F,r) + N(B(\alpha, R)\cap E, r)  
\]
and Lemma~\ref{equivforms1}.
\end{proof}

\begin{example}\label{irreducible}
Let $\alpha$ be an imaginary quadratic integer. Then the set of all irreducible elements of $\mathbb{Z}[\alpha]$ contains $(k,\epsilon, \mathbf{e})$-APs for every $k\geq 3$, $\epsilon>0$, and orthonormal basis $\mathbf{e}$ for $\mathbb{C}$.
\end{example} 
The following discussion is based on Seki's proof \cite{Seki}. He proved the sum of reciprocals of primes is divergent by using Roth's theorem \cite{Roth} which states that every subset of positive integers with positive upper density contains arithmetic progressions of length 3. This is a special case of Szmer\'edi's theorem \cite{Szemeredi}. 
 
\begin{proof}

By definition, we may  write $\mathbb{Z}[\alpha]=\{a+b\alpha \colon a,b\in \mathbb{Z}\}$. Let $P$ be the set of all irreducible elements of $\mathbb{Z}[\alpha]$. We may write $P=\{p_j : j=1,2,\cdots \}$ such that $|p_1|\leq|p_2|\leq\cdots$. It is enough to show that
\[
\sum_{j=1}^\infty \frac{1}{|p_j|^2}=\infty
\]
from Corollary~\ref{main5}. Let $E(R)= |B(0,R)\cap\mathbb{Z}[\alpha]|$. Then there exist constants $C_\alpha$ and $D_\alpha$ such that
\[
C_\alpha R^2 \leq E(R) \leq D_\alpha R^2.
\]
  Assume that $\sum_{j=1}^\infty \frac{1}{|p_j|^2}<\infty$. Then there exists $k\in \mathbb{N}$ such that 
\[
	\sum_{j=k+1}^\infty \frac{1}{|p_j|^2}\leq \frac{C_\alpha}{2D_\alpha}.
\]
Here let 
\begin{align*}
S&= \{x\in \mathbb{Z}[\alpha]\mid \text{ there exist }\alpha_1,\ldots,\alpha_k\geq0 \text{ such that }x=p_1^{\alpha_1}\cdots p_k^{\alpha_k} \}, \\
S(R)&= |B(0,R)\cap S|,\ S^c(R)= |B(0,R)\cap ( \mathbb{Z}[\alpha]\setminus S)|.
\end{align*}
Then we have
\begin{align*}
B(0,R)\cap ( \mathbb{Z}[\alpha]\setminus S)
&= \{x\in \mathbb{Z}[\alpha] \colon x \notin S,|x|\leq R \} \\
 & \subseteq \bigcup_{j=k+1}^\infty\big\{xp_j\in \mathbb{Z}[\alpha]\;:\; |x|\leq R/|p_j|\;,\;x\in\mathbb{Z}[\alpha] \big\}.
\end{align*}
This yields that
\[
S^c(R)\leq \sum_{j=k+1}^\infty E(R/|p_j|)\leq \sum_{j=k+1}^\infty D_\alpha R^2/|p_j|^2< C_\alpha R^2/2.
\]
Therefore we have $S(R)=E(R)-S^c(R)\geq C_\alpha R^2-C_\alpha R^2/2=C_\alpha R^2/2$, which implies that
\begin{equation}\label{example-f1}
	\Uzeta S=2
\end{equation}
from Remark~\ref{r-ind}. For every $\delta_1,\ldots,\delta_k \in \{0,1\}$ let 
\[
S(\delta_1,\ldots,\delta_k)=\{x\in S\;:\; \delta_1\equiv\alpha_1,\ldots,\delta_k\equiv \alpha_k \text{ (mod $2$)}\text{ and } x=p_1^{\alpha_1}\cdots p_k^{\alpha_k} \}.
\]
We see that 
\[
\bigcup_{\delta_1,\ldots,\delta_k \in \{0,1\}}S(\delta_1,\ldots,\delta_k)=S
\]
and 
\[
\frac{S(\delta_1,\ldots,\delta_k)}{p_1^{\delta_1}\cdots p_k^{\delta_k}}\subseteq \{x^2 : x\in \mathbb{Z}[i]\}
\]
where we define $\lambda F =\{\lambda x \colon x\in F \}$ for every $\lambda \in \mathbb{C}$ and $F\subseteq \mathbb{C}$. By (ii) in Lemma~\ref{union} and (\ref{example-f1}), we have 
\[
2+ \Uzeta S 
=\max_{\delta_1,\ldots,\delta_k \in \{0,1\}}\Uzeta S(\delta_1,\ldots,\delta_k) 
\]
 By the fact that the function $x \mapsto \lambda x$ is quasi-isometric and Theorem~\ref{quasi-iso}, we have
\[
\max_{\delta_1,\ldots,\delta_k \in \{0,1\}}\Uzeta S(\delta_1,\ldots,\delta_k)=\max_{\delta_1,\ldots,\delta_k \in \{0,1\}}\Uzeta \frac{S(\delta_1,\ldots,\delta_k)}{p_1^{\delta_1}\cdots p_k^{\delta_k}} 
\]
 Since we have
\[
\frac{S(\delta_1,\ldots,\delta_k)}{p_1^{\delta_1}\cdots p_k^{\delta_k}} \subseteq \{x^2 : x\in \mathbb{Z}[\alpha]\}
\]
for every $\delta_1,\ldots,\delta_k \in \{0,1\}$, we get
\[
\max_{\delta_1,\ldots,\delta_k \in \{0,1\}}\Uzeta \frac{S(\delta_1,\ldots,\delta_k)}{p_1^{\delta_1}\cdots p_k^{\delta_k}} \leq \Uzeta \{x^2 : x\in \mathbb{Z}[\alpha]\} 
\]
by (i) in Lemma~\ref{union}. Let $\Gamma = \{x^2 \colon x\in \mathbb{Z}[\alpha] \}$. By definition, $\zeta(\sigma ;0,\Gamma)=\zeta(2\sigma; 0, \mathbb{Z}[\alpha])$ holds. Thus we get  
\[
\Uzeta \{x^2 : x\in \mathbb{Z}[\alpha]\}=
 \frac{1}{2}\: \Uzeta \mathbb{Z}[\alpha]\leq \frac{1}{2}\Assouad \mathbb{C}=1. 
\]
by Theorem~\ref{main4} and Proposition~\ref{main1}. Hence 
we have $2\leq \Uzeta S \leq 1$. This is a contradiction.
\end{proof}

Choose $\alpha=i$ in Example~\ref{irreducible}. Then we can conclude that $\sum_{p\equiv 1 \text{ mod4 }}\frac{1}{p}=\infty\;$ by the above discussion. In fact, we define
\begin{itemize}
\item[(i)] $P_1=\{1-i,\ 1+i,\ -1-i,\ -1+i\}$; 
\item[(ii)] $P_2=\{p\in \mathbb{Z}\; \colon p \text{ is a rational prime of the form }4k+3\}$;
\item[(iii)] $P_3=\{q\in \mathbb{Z}[i]\: \colon |q|^2 \text{ is a rational prime of the form }4k+1\}$.
\end{itemize}
It is known that $P=P_1\cup P_2\cup P_3$ (see \cite[Theorem 252]{HardyWright}), where $P$ is the set of all irreducible elements of $\mathbb{Z}[i]$. It is clear that $\sum_{q\in P_1}|q|^{-2}<\infty$ and $\sum_{q\in P_2}|q|^{-2}<\infty$. Thus we have $\sum_{q\in P_3}|q|^{-2}=\infty$ by the proof of Example~\ref{irreducible}. Here let $r(p)=\#\{(x,y)\in \mathbb{Z}\times\mathbb{Z} \;:\;x^2+y^2=p \}.$ It is  well-known that $r(p)=8$ holds for every prime number $p$ such that   \cite[Theorem 251]{HardyWright}. Therefore we have for all $1<\sigma$
\[
\sum_{q\in P_3}\frac{1}{|q|^{2\sigma}}=\sum_{p\equiv 1 \text{ mod} 4}\frac{r(p)}{p^\sigma}=8\sum_{p\equiv 1 \text{ mod} 4}\frac{1}{p^\sigma}<8\sum_{p\equiv 1 \text{ mod}4}\frac{1}{p},
\] 
which implies $\sum_{p\equiv 1 \text{ mod4 }}\frac{1}{p}=\infty$
.

\begin{example}
Let $p(n)$ be the $n$-th prime number. The set of all primes of the form $p(p(n))$ contains $(k,\epsilon, \{1\})$-APs for every $k\geq 3$ and $\epsilon>0$.
\end{example}

\begin{proof}
Let $P=\{p(n)\colon n=1,2,\ldots\}$ and $Q=\{p(p(n))\colon n=1,2,\ldots \}$. It is enough to show that $\Uzeta Q=1$ from Corollary~\ref{main5}. Let $\pi(N)=|[1,N]\cap P|$. Here we have
\[
|Q\cap [1,N]|= |P\cap [1,\pi(N)]| = \pi(\pi(N)). 
\]
Thus by the prime number theorem, we obtain
\begin{equation}\label{distribution}
 \frac{|Q\cap [1,N]|(\log N)^2}{N} \longrightarrow 1 \quad \text{ as } N\longrightarrow \infty,
\end{equation}
which implies that $\Uzeta Q=1$ by (\ref{Euclidsp}).
\end{proof}

An element of the set $Q$ is called a super prime. Note that the sum of reciprocals of all super primes is convergent by (\ref{distribution}) and Abel's identity \cite[Theorem~4.2]{Apostol} but 
$
\sum_{q\in Q} q^{-\sigma}
$ 
is divergent for all $0<\sigma<1$.

\appendix

\section{Standard properties for zeta dimensions}\label{appendixA}

In this appendix, we show the existence of $r$-nets and standard properties of fractal dimensions for zeta dimensions. We find that the upper and lower zeta functions share some properties with upper and lower box dimensions, respectively.

\begin{proposition}
Let $(X,d)$ be a non-empty metric space For every $r>0$ there exists $M\subseteq X$ such that $M$ is an $r$-net.  
\end{proposition}

\begin{proof}
We define a family of sets ${\mathcal U} \subseteq \mathcal{P}(X)$ by
\[
	{\mathcal U}=\{A\subseteq X \mid A \text{ is $r$-discrete}\} \;.
\]
Since $X$ is non-empty, ${\mathcal U}$ is non-empty.  Here we take a subset ${\mathcal V}\subseteq {\mathcal U}$ which is totally ordered with respect to the inclusion. Let $C$ be the union of ${\mathcal V}$  i.e.
\[
	C=\bigcup_{A\in{\mathcal V}} A\;.
\]
Then $C$ is $r$-discrete. Indeed, for all distinct elements, $x,y\in C$ there exist $A_1,A_2 \in {\mathcal V}$ such that $x\in A_1$ and $y\in A_2$. We have either $A_1\subseteq A_2$ or $A_2\subseteq A_1$ since ${\mathcal V}$ is totally ordered. Without loss of generality, we may assume that $A_1\subseteq A_2$. Thus we have $x\in A_2$ and $y\in A_2$. Since $A_2$ is $r$-discrete, we have $d(x,y)\geq r$. Therefore $C$ is $r$-discrete.

It is clear that $C$ is an upper bound of ${\mathcal V}$. Therefore from Zorn's lemma, there exists a set $M\in {\mathcal U}$ such that $M$ is a maximal element in ${\mathcal U}$ with respect to the inclusion. Assume that there exists $x_0 \in X$ such that $d(x_0,y)\geq r$ for all $y\in M$. Then $M\cup\{x_0\}$ is clearly $r$-discrete, and we see that $M\subsetneq M\cup\{x_0\}$, which contradicts the maximality of $M$. Therefore $M$ is an $r$-net of $X$. 
\end{proof}

\begin{proposition}
Let $(X,d)$ be a non-empty metric space. Then the following properties hold:
\begin{itemize}
\item[(i)] For every $A\subseteq B \subseteq X$, we have
\[
\Uzeta A \leq \Uzeta B,\quad \Lzeta A\leq \Lzeta B. 
\]
\item[(ii)] For every $A, B\subseteq X$, we have
\[
 \Uzeta (A\cup B)= \max \{\Uzeta(A),\ \Uzeta(B)\}.
\]
\item[(iii)] Let $(Y,\ell)$ be a non-empty metric space. We have 
\[
 \Uzeta X + \Lzeta Y \leq \Uzeta (X\times Y) \leq \Uzeta X+ \Uzeta Y
\]  
and
\[
 \Lzeta X + \Lzeta Y  \leq \Lzeta( X \times Y) \leq \Lzeta X +\Uzeta Y,
\]
where we define the metric between elements $z_1=(x_1,y_1)$ and $z_2=(x_2, y_2)$ in $X\times Y$ as $\eta(z_1,z_2)=\max\{d(x_1,x_2),\ell(y_1,y_2)\}$.
\item[(iv)] For every $m\geq 1$, we have
\[
\Uzeta X^m = m\ \Uzeta X,\quad \Lzeta X^m =m\ \Lzeta X.
\]

\item[(v)] If $(X,d)$ be an $m$-dimensional Banach space, then we have
\[
 \Uzeta X= \Lzeta X=m.
\]  
\end{itemize}

\end{proposition}
The upper and lower zeta dimensions can be replaced by the upper and lower box dimensions in (i) to (iv) (see \cite{Falconer}). Note that the upper and lower zeta dimensions of any bounded set are zero. On the other hand, the upper and lower box dimensions of any unbounded set are infinity.

\begin{proof}
We have already proved (i) and (ii) of this proposition in Lemma~\ref{union}. We now show that (iii) of this proposition. Fix $\alpha_1\in X$ and $\alpha_2 \in Y$. Let $\alpha=(\alpha_1, \alpha_2)$. Fix any $R>0$ and $r>0$. Then we have
\begin{equation}\label{product1}
N(B(\alpha, R)\cap (X\times Y),r) \leq N(B(\alpha_1, R)\cap X,r) \cdot N(B(\alpha_2,R)\cap Y , r),
\end{equation}
since if $\mathcal{U}$ and $\mathcal{V}$ are finite $r$-covers of $B(\alpha_1, R)\cap X$ and $B(\alpha_2, R)\cap Y$, then 
\[
\{U\times V\colon U\in \mathcal{U}, V\in \mathcal{V}\}
\]
 is also a finite $r$-cover of $B(\alpha, R)\cap (X\times Y)$. By Lemma~\ref{equivforms1}, we obtain that 
\[
\Uzeta (X\times Y) \leq \Uzeta X +\Uzeta Y
\]
and 
\[
\Lzeta (X\times Y) \leq \Lzeta X +\Uzeta Y.
\]
Let us show that 
\begin{equation}\label{product}
N(B(\alpha_1, R)\cap X,4r) \cdot N(B(\alpha_2,R)\cap Y, 4r)\leq N(B(\alpha,R) \cap (X\times Y),r).
\end{equation}
In fact, if the right hand side is infinity, then it is clear. Thus we may assume that the right hand side is finite. Let $\Gamma$ and $\Omega$ be $(2r)$-nets of $X$ and $Y$, respectively. Then $\Gamma\times \Omega$ is also a $(2r)$-net of $X\times Y$. Here there exists a quasi-invarient surjection from $X\times Y$ to $\Gamma\times \Omega$. Therefore by Theorem~\ref{quasi-iso}, we have 
\[
\Uzeta (X\times Y)=\Uzeta(\Gamma\times \Omega),\quad \Lzeta (X\times Y) =\Lzeta (\Gamma\times \Omega).
\]

Let $N=N(B(\alpha,R) \cap (\Gamma\times \Omega),r)$. We can take an $r$-cover $\mathcal{U}$  of $B(\alpha,R) \cap (\Gamma\times \Omega)$ such that the number of elements in $\mathcal{U}$ is equal to $N$. Each $U\in \mathcal{U}$ intersects at most one point of $B(\alpha,R) \cap (\Gamma\times \Omega)$ since the distances of all distinct two elements of $\Gamma\times \Omega$ are at least $r$. Thus this implies that
\begin{align*}  
&\, |B(\alpha_1, R)\cap \Gamma|\cdot|B(\alpha_2,R)\cap \Omega|  
= |B(\alpha, R)\cap (\Gamma\times \Omega)| \\ 
&\quad \leq |\mathcal{U}| 
=N(B(\alpha,R) \cap (\Gamma\times \Omega),r) 
\end{align*}
Here the families of sets
\[
\{B(x,2r)\colon x\in B(\alpha_1, R)\cap \Gamma\} \text{ and }
 \{B(x,2r)\colon x\in B(\alpha_2, R)\cap \Omega\}
\]
are $(4r)$-covers of $B(\alpha_1, R)\cap X$ and $B(\alpha_2,R)\cap Y$, respectively. This yields that 
\[
N(B(\alpha_1, R)\cap X,4r) \cdot N(B(\alpha_2,R)\cap Y , 4r) \leq   |B(\alpha_1, R)\cap \Gamma|\cdot|B(\alpha_2,R)\cap \Omega|
\]
Therefore we obtain (\ref{product}). Hence by Lemma~\ref{equivforms1} and (\ref{product}), we have
\[
\Uzeta X +\Lzeta Y \leq \Uzeta(X\times Y)
\]
and 
\[
 \Lzeta X +\Lzeta Y\leq \Lzeta(X\times Y).
\]
Let us show (iv). By combining (\ref{product1}) and (\ref{product}), we see that
\[
N(B(\alpha_1,R)\cap X, 4^m r)^m \leq N(B(\alpha, R)\cap X^m,r) \leq N(B(\alpha_1, R)\cap X, r)^m
\]
for all $r>0$, where let $\alpha=(\alpha_1,\ldots, \alpha_1)\in X^m$. Hence by Lemma~\ref{equivforms1}, we obtain (iv). At last, let us show (v). For every $m$-dimensional Banach space $X$ there exists a quasi-isometric surjection from $X$ to $\mathbb{R}^m$. Therefore it is enough to show that 
\[
     \Uzeta \mathbb{R}=\Lzeta \mathbb{R} =1 
\]
from Theorem~\ref{quasi-iso} and (iv) in this proposition. 
Fix any $r>0$. It is obvious that 
\[
N(B(0,R)\cap \mathbb{R}, r)\geq \frac{1}{2r} R 
\]
for sufficiently large $R>0$. Therefore we have
\[
1\leq \Lzeta \mathbb{R} \leq \Uzeta \mathbb{R} \leq \Assouad \mathbb{R} =1
\]
by Lemma~\ref{equivforms1} and Theorem~\ref{main4}

\end{proof}

\section{Asymptotic dimension and zeta dimensions}\label{AsympZeta}

In this section, we show that the asymptotic dimension is not always an upper or lower bound for the zeta dimensions. More precisely, we will construct the following examples:
\begin{example}\label{exam1}
Let $A_n=\{j+2^n \colon j=0,1,\ldots, n-1\}$ for every $n=1,2,\ldots,$ and define $A=\bigcup_{n=1}^{\infty} A_n$. Then we have
\[
\asdim A=1,\quad \Lzeta A=\Uzeta A=0.
\]
\end{example}
\begin{example}\label{exam2}
Let $B(\alpha)=\{n^\alpha \colon n=1,2,\cdots \}$ for every $0<\alpha<1$. Then we have
\[
\asdim B(\alpha)=0,\quad \Lzeta B(\alpha)=\Uzeta B(\alpha)=1/\alpha.
\]
\end{example}
Let $(X,d)$ be a metric space. For every $r>0$, a family $\mathcal{U}$ of subsets of $X$ is called \textit{$r$-disjoint} if for all $U, V\in \mathcal{U}$ with $U\neq V$, we have $\mathrm{dist}(U,V)\geq r$, where we define 
\[
\mathrm{dist}(U,V)=\inf\{d(x,y)\colon x\in U, y\in V\}.
\]
Here we say that the \textit{asymptotic dimension} of a metric space $X$ does not exceed $n$ and write $\mathrm{asdim} X\leq n$ if for every $r<\infty$ there exist $r$-disjoint families $\mathcal{U}^0,...,\mathcal{U}^n$ of uniformly bounded subsets of $X$ such that $\bigcup_{i=0}^n\mathcal{U}^i$ is a cover of $X$. We say that the asymptotic dimension of $X$ is $n$ and write $\asdim X=n$ if $\asdim X \nleq n-1$ and $\asdim X\leq n$. By definition, it is clear that $\asdim F\leq \asdim E$ for every $E\subseteq F$. The asymptotic dimension is introduced by Gromov \cite{Gromov}. We can see several equivalent definitions of the asymptotic dimension in \cite{BellDranishnikov}. From Example~\ref{exam1} and Example~\ref{exam2}, neither $\asdim X \leq \Uzeta X$ nor $\Lzeta X \leq \asdim X$  holds for any metric space $X$. We do not find any relations between the asymptotic dimension and the zeta dimensions in this paper.

\begin{proof}[Proof of Example~\ref{exam1}]
Note that $A_1,A_2,\ldots$ are pairwise disjoint since
\[
 2^{m}+m-1<2^{n}
\] 
holds for all $1\leq m<n$. We firstly prove that $\Uzeta A=0$. Fix any positive integer $N$ and choose a positive integer $n$ such that $2^n \leq N <2^{n+1}$. Then we have
\[
|A\cap [1,N]|\leq |A_1|+\cdots+|A_n|
=1+\cdots+n\leq 2n^2\leq \frac{2}{(\log 2)^2} (\log N)^2,
\]
which implies that
\[
0\leq \Uzeta A =\limsup_{N\rightarrow \infty} \frac{\log |A\cap [1,N]|}{\log N }\leq  \limsup_{N\rightarrow \infty} \frac{\log \left(\frac{2}{(\log 2)^2} (\log N)^2 \right) }{\log N }= 0.
\]
We next show that $\asdim A=1$. It is enough to show that $\asdim A>0$ since
\[
\asdim A\leq \asdim \mathbb{R}=1.
\]
Assume that $\asdim A=0$. Then there exists a $2$-disjoint family $\mathcal{U}$ of uniformly bounded subsets of $A$ such that $A\subseteq \bigcup \mathcal{U}$. Choose a large positive integer $M$ such that the diameter of $U$ is at most $M$ for every $U\in \mathcal{U}$. Then let $a_0=2^{2M} \in A_{2M}$, and take a set $U\in \mathcal{U}$ such that $a_0\in U$. By the definition
 of $M$, there exists $j_0=1,2,\ldots, 2M-1$ such that $2^{2M}+j_0-1 \in U$ and $2^{2M}+j_0 \notin U$. Hence we can choose $V\in \mathcal{U}$ with $V \neq U$ such that $2^{2M}+j_0\in V$, which implies that 
\[
\mathrm{dist}(U,V)\leq 1<2. 
\]
This is a contradiction to the condition that $\mathcal {U}$ is a $2$-disjoint family. Therefore we have
\[
\asdim A>0.
\]
\end{proof}

\begin{proof}[Proof of Example~\ref{exam2}]
Fix any real number $\alpha\in (1,\infty)$. It is clear that $\Lzeta B(\alpha)=\Uzeta B(\alpha)=1/\alpha$ since 
\[
|B(\alpha)\cap [1,N]|=\lfloor N^{1/\alpha} \rfloor   
\]
holds for all $N\geq 1$. We show that $\asdim B(\alpha)=0$. Fix any real number $r>0$. Let $t=\lceil (r/\alpha)^{1/(\alpha-1)} \rceil $, $U_0=\{n^\alpha \colon n=1,2,\ldots ,t \}$ and $U_j=\{(t+j)^\alpha \}$ for every $j=1,2,\ldots$. Then $\{U_j\}_{j=0}^\infty$ is a family of uniformly bounded subsets of $B(\alpha)$, and 
\[
B(\alpha) \subseteq \bigcup_{j=0}^\infty U_j.
\]
Furthermore, by the mean value theorem, we find that
\[
\mathrm{dist}(U_j, U_{j+1}) = (t+j+1)^\alpha -(t+j)^\alpha \geq \alpha t^{\alpha-1}\geq r    
\]
for every $j=0,1,\ldots$. Therefore the family $\{U_j\}_{j=0}^\infty$ is $r$-disjoint. This yields that $\asdim B(\alpha)=0$.
\end{proof}

\section*{Acknowledgement}
This work is supported by Grantin-Aid for JSPS Research Fellow (Grant Number: 19J20878). The author would like to  thank Professor Kohji Mastumoto and Doctor Shin-ichiro Seki for useful comments.

\end{document}